\documentclass{amsart}
\usepackage{amssymb}
\def\bs{\boldsymbol}
\def\bg{\bs{g}}
\def\bm{\bs{m}}
\def\bbs{\bs{s}}
\def\BD{\bs{D}}
\def\BV{\bs{V}}
\def\bu{$\bullet$\quad}
\def\Cal{\mathcal}
\def\casesep{, & \text{ if }}
\def\CC{\mathbb C}
\def\CO{\mathcal O}
\def\CS{\mathcal S}
\def\DD{\mathbb D}
\def\eps{\varepsilon}
\def\Ree{\operatorname{Re}}
\def\ord{\operatorname{ord}}
\def\PSH{\mathcal{PSH}}
\def\phi{\varphi}
\def\RR{\mathbb R}
\def\too{\longrightarrow}
\def\tuu{\longmapsto}
\def\ZZ{\mathbb Z}
\makeatletter
\def\th@mytheorem{%
  \let\thm@indent\noindent
  \thm@headfont{\bfseries}
    \itshape
}
\def\th@myremark{%
  \let\thm@indent\noindent
  \thm@headfont{\bfseries}
}
\tagsleft@false
\makeatother
\theoremstyle{mytheorem}
\newtheorem{Theorem}{Theorem}
\theoremstyle{myremark}
\newtheorem{Remark}[Theorem]{Remark}

\begin{document}
\title{A remark on the Sibony function}
\author[M.~Jarnicki]{Marek Jarnicki}
\address{Jagiellonian University, Faculty of Mathematics and Computer Science, Institute of Mathematics,
{\L}ojasiewicza 6, 30-348 Krak\'ow, Poland}
\email{Marek.Jarnicki@im.uj.edu.pl}

\author[P.~Pflug]{Peter Pflug}

\address{Carl von Ossietzky Universit\"at Oldenburg, Institut f\"ur Mathematik, Postfach 2503, D-26111 Oldenburg, Germany}
\email{Peter.Pflug@uni-oldenburg.de}

\thanks{The research was partially supported by the OPUS grant no. 2015/17/B/ST1/00996 that was financed by the National
Science Centre, Poland}

\begin{abstract}
We present an effective formula for the Sibony function for all elementary Reinhardt domains.
\end{abstract}



\maketitle

For a domain $G\subset\CC^n$ let
$$
\bbs_G(a,z):=\sup\{\sqrt{u(z)}: u\in\CS_G(a)\},\quad a,z\in G,
$$
where
$$
\CS_G(a):=\{u:G\too[0,1): \log u\in\PSH(G),\;u(a) = 0,\; u\in\Cal C^2(\{a\})\}
$$
and $u\in\Cal C^2(\{a\})$ means that $u$ is $\Cal C^2$ in a neighborhood of $a$ (cf.~\cite{JarPfl2013}, \S\;4.2).
The function $\bbs_G$ is called the \emph{Sibony function} for $G$.

For $n\geq2$ and $\alpha=(\alpha_1,\dots,\alpha_n)\in\RR^n\setminus\{0\}$ let
$$
\BD_{\alpha}^n=\BD_\alpha:=\{z\in\CC^n(\alpha): |z^\alpha|:=|z_1|^{\alpha_1}\cdots|z_n|^{\alpha_n}<1\},
$$
where $\CC^n(\alpha):=\{(z_1,\dots,z_n)\in\CC^n: \forall_{j\in\{1,\dots,n\}}: (\alpha_j<0 \Longrightarrow z_j\neq0)\}$.
The domain $\BD_\alpha$ is called the \emph{elementary Reinhardt domain}.

Our aim is to present an effective formula for $\bbs_{\BD_\alpha}$. Partial results were presented in \cite{JarPfl2013}.

\medskip

Notice that $\BD_\alpha$ is a pseudoconvex $n$-circled domain, i.e.~if $(z_1,\dots,z_n)\in\BD_\alpha$, then
$(e^{i\theta_1}z_1,\dots,e^{i\theta_n}z_n)\in\BD_\alpha$ for arbitrary $\theta_1,\dots,\theta_n\in\RR$.
We say that $\BD_\alpha$ is of \emph{rational type} (resp.~\emph{irrational type}) if $\alpha\in\RR\cdot\ZZ^n$
(resp.~$\alpha\not\in\RR\cdot\ZZ^n$). It is clear that $\BD_\alpha=\BD_{t\alpha}$, $t>0$.

It is well-known that the system $(\bbs_G)_G$ is \emph{holomorphically contractible},
i.e.~if $F:G\too D$ is holomorphic (where $D\subset\CC^m$ is a domain), then
$\bbs_D(F(a),F(z))\leq\bbs_G(a,z)$, $a,z\in G$. In particular, the function $\bbs_{\BD_{\alpha}}$ is $n$-rotation invariant, i.e.
$$
\bbs_{\BD_\alpha}(a,z)=\bbs_{\BD_\alpha}((e^{i\theta_1}a_1,\dots,e^{i\theta_n}a_n), (e^{i\theta_1}z_1,\dots,e^{i\theta_n}z_n)),\quad
\theta_1,\dots,\theta_n\in\RR.
$$
Moreover,  $\bm_G\leq\bbs_G\leq\bg_G$, where
\begin{align*}
\bm_G(a,z):&=\sup\{|f(z)|: f\in\CO (G,\DD),\;f(a) =0\}, \\
\bg_G(a,z):&=\sup\Big\{u:G\too [0,1): \log u\in\PSH(G),\;\sup_{z\in G\setminus\{a\}}\frac{u(z)}{\|z-a\|}<+\infty\Big\};
\end{align*}
$\DD\subset\CC$ stands for the unit disc. The function $\bm_G$ (resp.~$\bg_G$) is called  the \emph{M\"obius pseudodistance}
(resp.~the \emph{pluricomplex Green function}). Both systems $(\bm_G)_G$ and $(\bg_G)_G$ are holomorphically contractible and
$\bm_{\DD}(a,z)=\bbs_{\DD}(a,z)=\bg_{\DD}(a,z)=\big|\frac{z-a}{1-\overline az}\big|$, $a,z\in\DD$.
Moreover,

\noindent (*) \quad if $\bg_G^2(a,\cdot)\in\Cal C^2(\{a\})$, then $\bbs_G(a,\cdot)=\bg_G(a,\cdot)$ (cf.~\cite{JarPfl2013}, Remark 4.2.8(b)).

The basic properties of $\bm_G$ and $\bg_G$ seem to be well understood. In contrast to that, almost nothing is known on $\bbs_G$.

In the case where $G=\BD_\alpha$ the following results are known.

\bu The rational case: If $\alpha_1,\dots,\alpha_n\in\ZZ$ are relatively prime, then
$$
\bm_{\BD_\alpha}(a,z)=\bm_{\DD}(a^\alpha,z^\alpha),\quad
\bg_{\BD_\alpha}(a,z)=\big(\bm_{\DD}(a^\alpha,z^\alpha)\big)^{1/r(a)},\quad a,z \in\BD_\alpha,
$$
where $r(a):=\begin{cases} 1 \casesep \sigma(a)=0 \\ \sum\limits_{\substack{j\in\{1,\dots,n\}:\\ \alpha_j>0,\;a_j=0}}
\alpha_j\casesep \sigma(a)\geq1\end{cases}$ (cf.~\cite{JarPfl2013}, Proposition 6.2.7),
$$
\sigma(n,\alpha,a)=\sigma(a):=\#\{j\in\{1,\dots,n\}: \alpha_j>0,\;a_j=0\}.
$$
Note that $r(a)=\ord_a(z^\alpha-a^\alpha)$.

One can easily prove that the function $\bg_{\BD_\alpha}^2(a,\cdot)\in\Cal C^2(\{a\})$ if and only if $\sigma(a)\leq1$.
Consequently, due to (*), if $\sigma(a)\leq1$, then
$$
\bbs_{\BD_\alpha}(a,z)=\big(\bm_{\DD}(a^\alpha,z^\alpha)\big)^{1/r(a)}=\bg_{\BD_\alpha}(a,z), \quad z\in\BD_\alpha.
$$
The most interesting case where $\sigma(a)\geq2$ remained unknown (since 1991).

\bu The irrational case:

$\bm_{\BD_\alpha}\equiv0$ (cf.~\cite{JarPfl2013}, Proposition 6.3.2);

$\bg_{\BD_\alpha}(a,z)=\begin{cases} 0 \casesep \sigma(a)=0 \\
|z^\alpha|^{1/r(a)} \casesep \sigma(a)=1 \end{cases}$ (cf.~\cite{JarPfl2013}, Proposition 6.3.3).
Hence, by (*), if $\sigma(a)\leq1$, then
$$
\bbs_{\BD_\alpha}(a,z)=\bg_{\BD_\alpha}(a,z),\quad z\in\BD_\alpha;
$$
cf.~\cite{JarPfl2013}, Proposition 6.3.10.
Once again, the case where $\sigma(a)\geq2$ remained unknown.

\begin{Theorem}\label{Thm1}
Assume that $\sigma(a)\geq2$. Let
$$
\mu(n,\alpha,a)=\mu(a):=\min\{\alpha_j: \alpha_j>0,\;a_j=0\}.
$$
Then
$$
\bbs_{\BD_\alpha}(a,z)=|z^\alpha|^{1/\mu(a)},\quad z\in\BD_\alpha.
$$
\end{Theorem}

\begin{Remark} As a consequence in the rational case we get:

\bu if $\sigma(a)\geq2$ and $\mu(a)=1$ (e.g.~$n=2$, $\alpha=(2,1)$, $a=0$), then
\begin{gather*}
\bm_{\BD_\alpha}(a,z)=\bbs_{\BD_\alpha}(a,z)<\bg_{\BD_\alpha}(a,z),\quad
z\in\BD_\alpha\setminus\BV_0,
\end{gather*}
where $\BV_0:=\{(z_1,\dots,z_n)\in\CC^n: z_1\cdots z_n=0\}$;

\bu if $\sigma(a)\geq2$ and $\mu(a)\geq2$ (e.g.~$n=2$, $\alpha=(3,2)$, $a=0$), then
\begin{gather*}
\bm_{\BD_\alpha}(a,z)<\bbs_{\BD_\alpha}(a,z)<\bg_{\BD_\alpha}(a,z),\quad
z\in\BD_\alpha\setminus\BV_0.
\end{gather*}

In the irrational case we have:

\bu if $\sigma(a)\geq2$, then
\begin{gather*}
0=\bm_{\BD_\alpha}(a,z)<\bbs_{\BD_\alpha}(a,z)<\bg_{\BD_\alpha}(a,z),\quad z\in\BD_\alpha\setminus\BV_0.
\end{gather*}
\end{Remark}

\begin{proof}[{Proof of Theorem \ref{Thm1}}]
Fix an $a=(a_1,\dots,a_n)\in\BD_\alpha$ with $\sigma(a)\geq2$.

If $\alpha_1\cdots\alpha_s\neq0$, $\alpha_{s+1}=\dots=\alpha_n=0$ for some $s\in\{1,\dots,n-1\}$, then
$\BD_\alpha^n=\BD_{\alpha'}^s\times\CC^{n-s}$ with $\alpha':=(\alpha_1,\dots,\alpha_s)$. Hence $\bbs_{\BD_\alpha^n}(a,z)=
\bbs_{\BD_{\alpha'}^s}(a',z')$ with $a':=(a'_1,\dots,a'_s)$, $z':=(z'_1,\dots,z'_s)$ (cf.~\cite{JarPfl2013}, Remark 4.2.9(b)).
Observe that $\sigma(n,\alpha,a)=\sigma(s,\alpha',a')$, $\mu(n,\alpha,a)=\mu(s,\alpha',a')$,
and $|z^\alpha|^{1/\mu(n,\alpha,a)}=|(z')^{\alpha'}|^{1/\mu(s,\alpha',a')}$.
\emph{This reduces the proof to the case where $\alpha_1\cdots\alpha_n\neq0$.}

If $\alpha_1,\dots,\alpha_s<0$ and $\alpha_{s+1},\dots,\alpha_n>0$ for some $s\in\{1,\dots,n\}$, then
consider the biholomorphic map
$$
(\CC\setminus\{0\})^s\times\CC^{n-s}\ni (z_1,\dots,z_n)\overset{F}
\tuu(1/z_1,\dots,1/z_s,z_{s+1},\dots,z_n)\in (\CC\setminus\{0\})^s\times\CC^{n-s}.
$$
Let $P:=\{(z_1,\dots,z_n)\in\CC^n: z_1\cdots z_s=0\}$, $\alpha':=(-\alpha_1,\dots,-\alpha_s,\alpha_{s+1},\dots,\alpha_n)$.
Note that $F$ maps biholomorphically $\BD_\alpha$ onto $\BD_{\alpha'}\setminus P$.
Hence
$$
\bbs_{\BD_\alpha}(a,z)=\bbs_{\BD_{\alpha'}\setminus P}(F(a),F(z))=\bbs_{\BD_{\alpha'}}(F(a),F(z)).
$$
Here we have used the following general property of $(\bbs_G)_G$: if $P\subset G$ is a closed pluripolar set, then
$\bbs_{G\setminus P}=\bbs_G|_{(G\setminus P)\times(G\setminus P)}$ (cf.~\cite{JarPfl2013}, Proposition 4.2.10(d)).
We have $\sigma(n,\alpha,a)=\sigma(n,\alpha',F(a))$, $\mu(n,\alpha,a)=\mu(n,\alpha',F(a))$,
and $|z^\alpha|^{1/\mu(n,\alpha,a)}=|(F(z))^{\alpha'}|^{1/\mu(n,\alpha',F(a))}$.
\emph{This reduces the proof to the case where $\alpha_1,\dots,\alpha_n>0$.}

Since the function $\bbs_{\BD_{\alpha}}$ is $n$-rotation invariant, we may assume that $a_1,\dots,a_n\geq0$.
We may also assume that $a_1,\dots,a_s>0$, $a_{s+1}=\dots=a_n=0$ with $s:=n-\sigma(a)$ and
$\mu(a)=\min\{\alpha_{s+1},\dots,\alpha_n\}=\alpha_n=1$. 

Observe that the function $z\tuu|z^\alpha|^2$ is of class $\Cal C^2(\{a\})$. Hence the function
$z\tuu|z^\alpha|$ belongs to $\CS_{\BD_\alpha}(a)$. Consequently,
$\bbs_{\BD_\alpha}(a,z)\geq|z^\alpha|$, $z\in\BD_\alpha$.

Let $\sqrt{u}$ be from $\CS_{\BD_\alpha}(a)$, $u\in\Cal C^2(\DD^n(a,\eps))$.
We have to prove that $\sqrt{u(z)}\leq|z^\alpha|$, $z\in\BD_\alpha$.
Since the set $\BV_0$ is analytic, it suffices to show that $\sqrt{u(z)}\leq|z^\alpha|$, $z\in\BD_\alpha\setminus\BV_0$
(cf.~e.g.~\cite{JarPfl2013}, B.4.23(c)). 
By the Liouville type theorem for plurisubharmonic functions we get $u=0$ on $\CC^{n-1}\times\{0\}$. Define
$$
\CC^{n-1}\times\DD\ni(\lambda_1,\dots,\lambda_n)\overset{F}\tuu (e^{\lambda_1},\dots,e^{\lambda_{n-1}},
\lambda_ne^{-(\alpha_1\lambda_1+\dots+\alpha_{n-1}\lambda_{n-1})})\in\CC^n.
$$
Observe that $|(F(\lambda))^\alpha|=|\lambda_n|$. In particular, $F(\CC^{n-1}\times\DD)\subset\BD_\alpha$.
Moreover, $F$ is surjective on $\BD_\alpha\cap(\CC_\ast^{n-1}\times\CC)$. In fact, given
$z=(z_1,\dots,z_n)\in\BD_\alpha\cap(\CC_\ast^{n-1}\times\CC)$ we first choose $\lambda_1,\dots,\lambda_{n-1}\in\CC$ such that
$z_j=e^{\lambda_j}$, $j=1,\dots,n-1$, and next we find a $\lambda_n\in\CC$ such that
$z_n=\lambda_ne^{-(\alpha_1\lambda_1+\dots+\alpha_{n-1}\lambda_{n-1})}$. Since $z\in\BD_\alpha$ we conclude that
$$
|\lambda_n|=|z_ne^{\alpha_1\lambda_1+\dots+\alpha_{n-1}\lambda_{n-1}}|
=|z_n|e^{(\alpha_1\Ree\lambda_1+\dots+\alpha_{n-1}\Ree\lambda_{n-1})}|=|z^\alpha|<1.
$$

Using once again the Liouville type theorem for plurisubharmonic functions we get
$\sqrt{u(F(\lambda_1,\dots,\lambda_n))}=v(\lambda_n)=\sqrt{u(F(0,\dots,0,\lambda_n))}=\sqrt{u(1,\dots,1,\lambda_n)}$.
Hence $v(0)=0$ (because $u=0$ on $\CC^{n-1}\times\{0\}$) and $v$ is log-subharmonic in $\DD$.

Fix $\lambda_1^0,\dots,\lambda_{n-1}^0\in\CC_\ast$ with $|e^{\lambda_j^0}-a_j|<\eps$, $j=1,\dots,n-1$. Then
$$
|\lambda_ne^{-(\alpha_1\lambda_1^0+\dots+\alpha_{n-1}\lambda_{n-1}^0)}|<\eps,\quad |\lambda_n|\ll1.
$$
Thus $v^2\in\Cal C^2(\{0\})$. The function $\DD\setminus\{0\}\ni\zeta\tuu\frac{v^2(\zeta)}{|\zeta|^2}$ is subharmonic and locally
bounded near $0$. Hence, it extends to a subharmonic function in $\DD$ and, by the maximum principle,  $v^2(\zeta)\leq|\zeta|^2$, $\zeta\in\DD$. Consequently, if
$z=F(\lambda)\in \BD_\alpha\cap(\CC_\ast^{n-1}\times\CC)$, then
$\sqrt{u(z)}=\sqrt{u(F(\lambda))}\leq|\lambda_n|=|z^\alpha|$.
\end{proof}

\bibliographystyle{amsplain}

\end{document}